\documentclass[a4paper]{article}            
\usepackage[all]{xy}\usepackage[latin1]{inputenc}        
\usepackage[dvips]{graphics}
\usepackage[dvips]{graphicx}
\usepackage{amsfonts}
\usepackage{amssymb}
\usepackage{amsmath}
\usepackage{amstext}
\usepackage{amsbsy}
\usepackage{amsopn}
\usepackage{amsthm}
\usepackage{amscd}
\usepackage{amsxtra}
\usepackage{mathrsfs}
\usepackage{upref}
\numberwithin{equation}{section}              
\newtheorem{theorem}{Theorem}[section]
\newtheorem{lemma}{Lemma}[section]
\newtheorem{proposition}{Proposition}[section]
\newtheorem*{proposition*}{Proposition}

\newtheorem*{corollary*}{Corollary}

\newtheorem*{definitions*}{Definitions}
\newtheorem*{conjecture*}{\bf Conjecture}

\newtheorem*{example*}{\bf Example}
\theoremstyle{remark}
\newtheorem{remark}{\bf Remark}[section]

\newcommand{\N}{\mathbb{N}}               

\newcommand{\R}{\mathbb{R}}
\newcommand{\ds}{\displaystyle} 
\mathchardef\varepsilon="010F
\mathchardef\epsilon="0122
\mathchardef\vartheta="0112
\mathchardef\varrho="011A
\mathchardef\rho="0125
\mathchardef\varphi="011E     
\mathchardef\phi="0127
\renewcommand 
\emptyset \varnothing
\begin{document}  
\date{}                                     
\title{On some nonlinear partial differential equations involving the $1-$Laplacian}
\author{\sc Mouna Kra\"iem \\ 
{\small Universit\'e de Cergy Pontoise , Departement de Mathematiques,}\\ 
{\small 2, avenue Adolphe Chauvin, 95302 Cergy Pontoise Cedex, France} \\ 
{\small e-mail: Mouna.Kraiem@math.u-cergy.fr}\\
{\small Acc\'ept\'e pour publication aux Annales de la Facult\'e des Sciences de Toulouse}}
\maketitle

\begin{abstract}

In this paper we present an approximation result concerning the first eigenvalue of the $1$-Laplacian operator. More precisely,  for $\Omega$ a bounded regular  open domain, we consider  a minimisation of the functional ${\ds \int_\Omega}|\nabla u|+n\left( {\ds \int _\Omega} |u|-1\right)^2
$ over the space $W_0^{1,1}(\Omega)$. For $n$ large enough, the infimum is achieved in some sense on $BV(\Omega)$, and letting $n$ go to infinity  this provides an approximation of the first eigenfunction for the first eigenvalue, since the term $n\left({\ds \int_\Omega} |u|^2-1\right)^2$ ``tends" to the constraint $\|u\|_1=1$. 
\end{abstract}

\section{Introduction: the first eigenvalue for the $1$-Laplacian}

In recent fields, several authors were interested on the study of the ``first eigenvalue" for the $1$-Laplacian operator, that we shall denote as the not everywhere defined $u \mapsto  -{\rm div} (\frac{\nabla u}{|\nabla u|})$.  

Due to the singularity of this operator, the definition of the first eigenvalue can  be  correctly  defined with the aid of  a variational  formulation: 
let $\lambda_1$ be defined as  
\begin{equation}\label{1Lap:eq:lambda_1}
\lambda_1:= \inf_{\begin{smallmatrix}
u\in W_0^{1,1}(\Omega)\\
\|u\|_1=1
 \end{smallmatrix}} \int_\Omega |\nabla u|.
\end{equation}
Notice that  $\lambda_1$ is well defined and is  positive, due to Poincar\'e's inequality. 

In order to justify the term  ``eigenvalue" for $\lambda_1$, 
one must prove the existence of an associated ``eigenfunction". 
As in the $p$-Laplacian case, an eigenfunction will be a solution of (\ref{1Lap:eq:lambda_1}).
Unfortunately, since $W^{1,1}(\Omega)$ is not a reflexif space, 
one cannot hope to obtain a solution for (\ref{1Lap:eq:lambda_1}) by classical arguments. 

This difficulty can be overcome by introducing the space $BV(\Omega)$, 
which is the weak closure of $W^{1,1}(\Omega)$, 
and by extending the infimum to that space, using the features of $BV(\Omega)$:  Density of regular maps in $BV$, existence of the trace map on the boundary... 
However, these properties  are  not sufficient to obtain solutions by classical methods, 
since the trace map --which is well defined on $BV(\Omega)$-- is not continuous  for the weak topology. 
This new difficulty can be ``solved" by introducing
--as it is the case in the theory of minimal surfaces and in plasticity and also for related problems-- 
a ``relaxed" formulation for (\ref{1Lap:eq:lambda_1}). 
This relaxed formulation consists in replacing the condition $\{u=0\}$ on the boundary by the addition of a term ${\ds \int_{\partial \Omega}} |u|$ in the functional to minimize. The new formulation is   then 

\begin{equation}\label{1Lap:eq:I_BV}
\inf_{\begin{smallmatrix}
u\in BV(\Omega)\\
\|u\|_1=1 \end{smallmatrix}}
 \int_\Omega |\nabla u|+\int_{\partial\Omega} |u|. 
\end{equation}
This problem has an infimum equal to $\lambda_1$. 
It can be seen by approximating function in $BV(\Omega)$ by functions in $W^{1,1}(\Omega)$ for a topology related to the narrow topology of measures.
This topology is precised in section 2.

Then  the existence of a minimizer of (\ref{1Lap:eq:I_BV}) in $BV(\Omega)$  can be proved, 
using classical arguments, arguments which  will be precised later in this paper.

To obtain the partial differential equation satisfied by a minimizer of (\ref{1Lap:eq:I_BV}), equation which can be seen as an eigenvalue's equation,  the author used in \cite{DF5} an approximation of (\ref{1Lap:eq:lambda_1}) by the following problem on $W_0^{1,1+\epsilon}(\Omega)$: 
\begin{equation}\label{1Lap:eq:lambda_epsilon}
\lambda_{1+ \epsilon} := \inf_
{\begin{smallmatrix}
u\in W_0^{1,1+\epsilon}(\Omega )\\ 
\|u\|_1=1 \end{smallmatrix}}
\int_\Omega|\nabla u|^{1+\epsilon},
\end{equation}
and proves that $\lambda_{1+\epsilon}$ converges to $\lambda_1$. 
Moreover, if $u_\epsilon$ is a positive solution of  the minimizing problem defined in (\ref{1Lap:eq:lambda_epsilon}), $u_\epsilon$ converges weakly in $BV(\Omega)$ to some $u$ which satisfies 
$$-{\rm div} \left(\frac{\nabla u}{|\nabla u|}\right) = \lambda_1,$$
in a sense which needs  to be precised, and is detailed in the present paper. 

Let us note that it is also proved in \cite{DF3} that there are caracteristic functions of  sets which are solutions. These sets are therefore called eigensets. 

Another approach is used in \cite{ACCh1}, \cite{ACCh2}, where the authors use the concept of Cheeger sets \cite{Cheeger}. 
In these papers, the authors present a remarkable construction of eigensets  for $2-$dimensional convex sets $\Omega$. 
Among their results, there is  the uniqueness of eigensets in the  case $N=2$. 

Our aim in the present article is to propose an approach of the first eigenvalue and the first eigenfunction of the $1$-Laplacian operator, using a penalization method, which consists in replacing the condition  ${\ds \int_{\Omega}}|u|=1$  in (\ref{1Lap:eq:I_BV}) by  the introduction of  the  term  $n\left({\ds \int_{\Omega}}|u|-1\right)^2$ inside the functionnal to minimize. 
This provides in the same time,  a new proof of the existence and uniqueness of a positive eigenfunction.
 
\section{Survey on known results about the space  $BV(\Omega)$}

We begin to  recall the definition of the space of functions with bounded variation. 
Let $\Omega$ be an open regular domain  in $\R^N, \ N>1$, and let ${\cal M}^1(\Omega)$ be the space of bounded measures in $\Omega$. 
We define 
$$BV(\Omega) = \left\{ u\in L^1(\Omega), \nabla u\in {\cal M}^1(\Omega)\right\}.$$
Endowed with the norm 
${\ds \int_\Omega} |\nabla u|+ {\ds \int_\Omega} |u|$, the space  $BV(\Omega)$ is a Banach space.

More useful is the weak  topology   for  variational technics :

We define {\it the weak topology} with  the aid of sequences,  as follows: we say that a sequence 
$u_n \rightharpoonup u$ weakly in $BV(\Omega)$  if the following two conditions are fullfilled:
\begin{itemize}
\item 
${\ds \int_\Omega} |u_n-u|\longrightarrow 0
 \quad \text{in} \ L^1(\Omega) \quad \text{when} \ n\longrightarrow\infty$,
 \item 
${\ds \int_\Omega}  \partial_i u_n \ \varphi \longrightarrow {\ds \int_\Omega} \partial_i u \ \varphi, 
\quad \forall i=1, 2,..., N \quad \forall \varphi \in {\cal C}_c(\Omega) \quad \text{when} \ n\longrightarrow \infty.$
\end{itemize} 
Let us note that the second convergence is also denoted as  {\it the vague convergence} of $\nabla u_n$ towards $\nabla u$. 

We shall also use the concept of {\it tight convergence} in $BV(\Omega)$: we say that a sequence 
$u_n$ converges tightly to $u$ in $BV(\Omega)$  if the following two conditions are fullfilled:
\begin{itemize}
\item 
$u_n \rightharpoonup u , \ \text{weakly} \quad \text{in} \ BV(\Omega)  \quad \text{when} \ n\longrightarrow \infty,$

\item 
${\ds \int_\Omega} |\nabla u_n| \longrightarrow {\ds \int_\Omega}  |\nabla u| \quad \text{when} \ n\longrightarrow\infty.$
\end{itemize}
Let us note that the last assertion is equivalent to say that, for all $\varphi\in {\cal C}(\bar\Omega, \R^N)$,
$$\int_\Omega \nabla u_n \cdot \varphi \longrightarrow \int_\Omega \nabla u \cdot \varphi,  \quad \text{when} \ n\longrightarrow\infty.$$

We now recall some facts about embedding and compact embedding from $BV(\Omega)$ to other $L^q$ spaces : 
\begin{itemize}
\item If $\Omega$ is  an open ${\cal C}^1$ set, then 
$BV(\Omega)$ is continuously embedded in $L^p(\Omega)$ for all $p\leq \frac{N}{N-1}$. 
\item If $\Omega$ is also bounded and smooth, the embedding is compact in $L^p(\Omega)$ for every $p < \frac{N}{N-1}.$
\end{itemize}

Finally we recall the existence of   a map, called {\it trace map},  defined on $BV(\Omega)$, which co\"incides with the restriction on $\partial \Omega$ of $u$
when $u$ belongs to ${\cal C}(\bar \Omega) \cap BV(\Omega)$ or less classically when $u\in W^{1,1}(\Omega)$. 
This map is continuous for  the strong topology,  and  is not continuous under the weak topology.  However the following property holds:
if $u_n\rightarrow u$ tightly in $BV(\Omega)$, then  
$$\int_{\partial \Omega} |u_n-u|\longrightarrow 0 \quad \text{for} \ n \rightarrow \infty.$$

We end this section by enouncing  a generalization of the Green's formula :  
this will allow us to give sense to the product $\sigma \cdot \nabla u$ when $\sigma$ is  in $ L^{\infty}(\Omega, \R^N)$, ${\rm div}\ \sigma \in L^N(\Omega)$ and $u\in BV(\Omega)$, and  
 will be useful  to give sense to  the partial differential equation associated to the eigenvalue.

Let us recall that  $\mathcal D(\Omega)$ is the space of ${\cal C}^\infty$-functions, with support on $\Omega$. 
\begin{proposition}\label{1Lap:prop:1.1}
Let  $\sigma \in L^{\infty}(\Omega, \R^N)$, ${\rm div}\sigma \in L^N(\Omega)$ and $u \in BV(\Omega)$.
Define the distribution $\sigma \cdot \nabla u$ by the following formula : for  any $\phi \in \mathcal D(\Omega)$,

\begin{equation}\label{1Lap:eq:1.5}
\langle \sigma \cdot \nabla u, \phi \rangle=-\int_\Omega ({\rm div}\sigma) u \phi-\int_\Omega (\sigma \cdot \nabla \phi) \ u .
\end{equation}
Then
$$| \langle \sigma \cdot \nabla u, \phi \rangle| \leq \|\sigma\|_\infty \langle| \nabla u|,| \phi| \rangle .$$
In particular, $ \sigma \cdot \nabla u$ is a bounded measure which satisfies:
$$|\sigma \cdot \nabla u|\leq \|\sigma\|_\infty |\nabla u|.$$
In addition, if $\phi \in \mathcal C(\overline{\Omega}) \cap  \mathcal C^1(\Omega)$, the following Green's Formula holds:
\begin{equation}\label{1Lap:eq:1.6}
\langle \sigma \cdot \nabla u, \phi \rangle
=-\int_\Omega ({\rm div} \sigma)  u \phi-\int_\Omega (\sigma \cdot \nabla \phi) \ u +\int_{\partial \Omega} \sigma \cdot \overrightarrow {n} u \ \phi,
\end{equation}
where $\overrightarrow {n}$ is the unit outer normal to $\partial \Omega$.
 
Suppose that $U \in BV(\R^N\setminus\overline{\Omega})$, that $u \in BV(\Omega)$ and define the function $\widetilde u$ as:
$$\widetilde u =
\begin{cases}
u \ {\rm in} \ \Omega,\\
U \ {\rm in} \ \R^N\setminus\overline{\Omega}.\\
\end{cases}$$
Then $\widetilde u \in BV(\R^N)$ and
$$\nabla \widetilde u = \nabla u \ \chi_{\Omega} +\nabla U \ \chi_{(\R^N\setminus \overline{\Omega})}+(U-u) \ \delta_{\partial \Omega},$$
where in the last term, $U$ and $u$ denote the trace of $U$ and $u$ on $\partial \Omega$ and $\delta_{\partial \Omega}$ denotes the uniform Dirac measure on  $\partial \Omega$.
Finally,  we introduce the measure $\sigma \cdot \nabla\widetilde u$ on $\overline{\Omega}$ by the formula
$$(\sigma \cdot \nabla \widetilde u)=(\sigma \cdot \nabla u)\chi_{\Omega} +\sigma \cdot\overrightarrow {n}(U-u) \ \delta_{\partial \Omega} $$
where $(\sigma \cdot \nabla u)\chi_{\Omega} $ has been defined in  \ref{1Lap:eq:1.5}. 
Then $\sigma \cdot \nabla \widetilde u$ is absolutely continuous with respect to $|\nabla \widetilde u|$, with the inequality
$$|\sigma \cdot \nabla \widetilde u| \leq \|\sigma\|_\infty |\nabla  \widetilde u|.$$
\end{proposition}
For a proof the reader can consult \cite{DF2}, \cite{Ko.Tem}, \cite{Str.Tem1}.

\section{Presentation of the  results}
We now describe the approximation result here enclosed. 
For  $n\in \N^*$, let us consider the following minimization problem:
\begin{equation}\label{1Lap:eq:lambda_n}
\lambda_{1, n}=\inf_{\begin{smallmatrix} 
u \in W^{1,1}_0(\Omega)  
\end{smallmatrix}}
\left\{\int_\Omega |\nabla u| +  n \left(\int_\Omega|u| -1 \right)^2\right\}.
\end{equation}
As it is done for analogous problem in \cite{}, let us introduce the relaxed formulation  associated : 
\begin{equation}\label{1Lap:eq:lambda_BV}
\widetilde \lambda_{1, n}= \inf_{u\in BV(\Omega)}\left\{\int_\Omega |\nabla u| +\int_{\partial \Omega} |u|+   n \left(\int_\Omega|u| -1 \right)^2\right\}. 
\end{equation}
We shall  prove in the following section the  result : 

\begin{theorem}\label{1Lap:thm:1.1}
Let $\Omega$ be a piecewise $\mathcal C^1$ bounded domain in $\R^N, \ N>1$.
For every $n\in \N^*$, the problem {\rm(\ref{1Lap:eq:lambda_BV})} possesses a  solution $u_n$ in $BV(\Omega)$ which can be chozen nonnegative.  
Moreover, $u_n$ satisfies the following partial differential equation: 
\begin{equation}\label{1Lap:eq:P_n}
\begin{cases}
-{\rm div} \ \sigma_n +2n \left({\ds \int}_ \Omega u_n -1\right) {\rm sign}^+ (u_n) =0\quad \text{\rm in} \ \Omega, \\
\sigma_n \in  L^{\infty}( \Omega,\R^N), \ \|\sigma_n\|_{\infty} \leq 1,\\
\sigma_n \cdot \nabla u_n= |\nabla u_n|  \quad \text{\rm in} \ \Omega, \\
u_n \  \text{\rm is not identically zero}, -\sigma_n \cdot \overrightarrow {n} (u_n)= u_n \quad \text{\rm on}\ \partial\Omega,
\end{cases}
\end{equation}
where $\overrightarrow {n}$ denotes the unit outer normal to $\partial\Omega$,  $\sigma_n \cdot\nabla u_n$ is the measure defined in Proposition \ref{1Lap:prop:1.1} and ${\rm sign}^+ (u_n)$ is some function in $L^\infty(\Omega)$ such that ${\rm sign}^+ (u_n) u_n=u_n \ \text{in} \ \Omega$.

Moreover $\lambda_{1, n}$ converges towards  $\lambda_1$ and $u_n$ converges towards the first eigenfunction $u$.
\end{theorem}

\begin{remark}
Clearly, $u_n$ is not identically zero for $n$ large enough as soon as $n>\lambda_1$.
\end{remark}

\begin{remark}
From Proposition \ref{1Lap:prop:1.1} (with $U=0$), the conditions
$$ \sigma_n \cdot \nabla u_n= |\nabla u_n|  \quad \text{\rm in} \ \Omega \ , \  \quad \ -\sigma_n \cdot \overrightarrow {n} (u_n)= u_n \quad \text{\rm on}\ \partial\Omega,$$ 
are equivalent to
$$\sigma_n \cdot \nabla \widetilde u_n =|\nabla \widetilde u_n| \quad \text{\rm on}\ \Omega \cup \partial \Omega.$$
\end{remark}

\begin{remark}
The identity $ \sigma_n \cdot \nabla u_n= |\nabla u_n|$ makes sense since
$$-{\rm div} \ \sigma_n= -2n \left(\int_ \Omega u_n -1\right) {\rm sign}^+(u_n),$$
which implies that  ${\rm \ div} \sigma_n \in L^{\infty}(\Omega)$, therefore $\sigma_n \cdot \nabla u_n$ is well-defined by Proposition \ref{1Lap:prop:1.1}.
\end{remark}

\bigskip

We subdivide the proof of Theorem \ref{1Lap:thm:1.1} into several steps : 

\begin{itemize}

\item First step: 
We use some  kind of regularization of the minimization problem by introducing    for  some $\epsilon>0$ and small 
$$\inf_{\begin{smallmatrix} 
u \in W^{1,1+\epsilon}_0(\Omega)
\end{smallmatrix}} 
\left\{\int_\Omega |\nabla u |^{1+\epsilon} +  n \left(\int_\Omega|u|^{1+\epsilon} -1 \right)^2 \right \}.$$
We prove that for $n$ large enough, this problem possesses a  solution which can be chozen nonnegative  and denoted by  $u_{n, \epsilon}$, which satisfies 
$$
\begin{cases}
-{\rm div} (|\nabla u_{n,\epsilon}|^{\epsilon-1} \nabla u_{n,\epsilon}) +2n \left({\ds \int}_ \Omega u_{n,\epsilon}^{1+\epsilon} -1 \right) u_{n,\epsilon}^\epsilon =0,  \quad \text{in} \ \Omega,\\
\end{cases}$$

\item Second step: 
 We extend $u_{n,\epsilon}$ by zero outside of $\Omega$ and observe that  the sequence still denoted  $(u_{n, \epsilon})$ is uniformly bounded in $BV(\R^N)$, more precisely 
$$\int_{\R^N} |\nabla u_{n,\epsilon}|^{1+\epsilon} \leq C.$$
Then we can extract from $u_{n,\epsilon}$ a subsequence, such that 
$u_{n,\epsilon}\rightharpoonup u_n$ weakly in $BV(\R^N)$. The limit function belongs to $BV(\R^N)$ and is zero outside of $\bar\Omega$. 

\item Third  step:
we prove that $\sigma_{n, \epsilon}=|\nabla u_{n,\epsilon}|^{\epsilon-1} \nabla u_{n,\epsilon}$ is uniformly bounded in $L^q(\Omega) \ \forall \ q<\infty$. Then we can extract from $\sigma_{n,\epsilon}$ a subsequence, such that 
$\sigma_{n,\epsilon}\rightharpoonup \sigma_n$ weakly in $L^q(\Omega) \ \forall \ q<\infty$, such that $\|\sigma\|_\infty\leq 1$ and $\sigma_n \cdot \nabla u_n= |\nabla u_n|$ in $\Omega \cup \partial\Omega$.

\item Fourth step:
we prove that $u_n$ is a solution of the minimizing problems (\ref{1Lap:eq:lambda_BV}) and (\ref{1Lap:eq:P_n}).
We also prove that $\sigma_n$ satisfies the problem (\ref{1Lap:eq:P_n}).

\item Fifth step:
we  establish that $\lambda_{1, n}$ converges strongly to $\lambda_1$ when $n$ goes to $\infty$ and that $u_n$ converges strongly to the first eingenfunction associated to $\lambda_1$.
\end{itemize}

\section{Proof of the main result}

We provide here the proof of Theorem \ref{1Lap:thm:1.1}, outlined as above.
\paragraph{Step 1:}
We prove here the existence and uniqueness of a positive solution for the following approximation problem 
\begin{equation}\label{1Lap:eq:2.2}
\lambda_{1+\epsilon,  n}=\inf_{\begin{smallmatrix} u \in 
W^{1,1+\epsilon}_0(\Omega)
\end{smallmatrix}}I_{1+\epsilon, n}(u), 
\end{equation} 
where $I_{1+\epsilon,  n}$ is the following functional
\begin{equation}\label{1Lap:eq:2.3}
I_{1+\epsilon, n}(u) =\int_\Omega |\nabla u |^{1+\epsilon} +  
n \left(\int_\Omega|u|^{1+\epsilon} -1 \right)^2,
\end{equation} 
for some positive $\epsilon$ given.

We first prove that $\lambda_{1+\epsilon,  n}$ is achieved, using standard variational technics:
Let $(u_i)_i$ be a minimizing sequence for $\lambda_{1+\epsilon,  n}$. 
Without loss of generality, up to replace $u_i$ by $|u_i|$, 
one may assume that $u_i$ is nonnegative.
Since $I_{1+\epsilon, n}$ is coercive, 
$(u_i)$ is  bounded in $W^{1,1+\epsilon}_0 (\Omega)$.

As a consequence, we may extract from it a subsequence, 
still denoted $(u_i)_i$, which converges weakly in $W^{1,1+\epsilon}_0 (\Omega)$ to some function $u_{n,\epsilon} \in W^{1,1+\epsilon}_0 (\Omega)$.
Furthermore, by the Rellich-Kondrakov Theorem \cite{Aub2}, \cite{Aub1}, \cite{Adams}, $(u_i)_i$ converges to $u_{n,\epsilon}$ in $L^{1+\epsilon}(\Omega)$.

Using the weak lower semicontinuity of the semi-norm $\int_\Omega |\nabla u|^{1+\epsilon}$  for the weak topology of $W^{1,1+\epsilon}_0 (\Omega)$, one has:
\begin{align*}
\lambda_{1+\epsilon,  n} 
&\leq \int_\Omega |\nabla u_{n,\epsilon} |^{1+\epsilon} +  
n \left( \int_\Omega|u_{n,\epsilon}|^{1+\epsilon} -1  \right)^2 \\ 
&\leq \liminf_{i\rightarrow+\infty} \left[\int_\Omega |\nabla u_i|^{1+\epsilon} + 
n\left(\int_\Omega|u_i|^ {1+\epsilon}-1 \right)^2\right]= 
\lambda_{1+\epsilon,  n}.
\end{align*}
Hence, $u_{n, \epsilon}$ is a solution of the minimization problem (\ref{1Lap:eq:2.2}).

We now prove that this weak solution solves the following partial differential equation:
\begin{equation}\label{1Lap:eq:P_{n,epsilon}}
\begin{cases}
-{\rm div} \sigma_{n, \epsilon}+
2n \left({\ds \int}_ \Omega u_{n,\epsilon}^{1+\epsilon} -1 \right)  
u_{n,\epsilon}^\epsilon=0 \quad \text{in} \ \Omega, \\
\sigma_{n,\epsilon}\cdot \nabla u_{n,\epsilon} =
|\nabla u_{n,\epsilon}|^{1+\epsilon}
\quad \text{in} \ \Omega, \\
 u_{n,\epsilon} > 0\quad \text{\rm in} \ \Omega, \ u_{n,\epsilon}=0 
\quad \text{\rm on} \ \partial\Omega.
\end{cases}
\end{equation}
Indeed, for every $h \in \mathcal D(\Omega)$, we have:
\begin{align*}
&D I_{1+\epsilon, n}(u_{n, \epsilon})\cdot h\\
&= (1+\epsilon) \left[\int_\Omega |\nabla u_{n,\epsilon} |^{\epsilon-1}\nabla u_{n,\epsilon} \cdot \nabla h +  
 2n \left(\int_ \Omega u_{n,\epsilon}^{1+\epsilon} -1 \right) 
 \int_ \Omega u_{n, \epsilon}^\epsilon h \right]\\
&=(1+\epsilon){\ds \int_\Omega}\left[-{\rm div} \left(|\nabla u_{n,\epsilon}|^{\epsilon-1} \nabla u_{n,\epsilon}\right) +2n \left(\int_ \Omega u_{n,\epsilon}^{1+\epsilon }-1 \right) u_{n,\epsilon}^\epsilon 
\right] h\\
 &=0 .
\end{align*}
Thus, we get:
\begin{equation}\label{1Lap:eq:2.5}
-{\rm div} \left(|\nabla u_{n,\epsilon}|^{\epsilon-1} \nabla u_{n,\epsilon}\right) +2n \left(\int_ \Omega u_{n,\epsilon}^{1+\epsilon }-1 \right) u_{n,\epsilon}^\epsilon =0, 
\end{equation}
in a distribution sense.

Since $u_{n,\epsilon}$ is a weak solution of equation (\ref{1Lap:eq:2.5}), by regularity results (as developped by Guedda-Veron \cite{GV}, see also Tolksdorf \cite{Tolk}), one gets that  $u_{n,\epsilon} \in   \mathcal C^{1,\alpha}(\bar \Omega) \ \forall \ \alpha \in (0,1)$. 
Moreover, since $u_{n, \epsilon}$ is a nonnegative weak solution of the equation (\ref{1Lap:eq:2.5}), by the strict maximum principle of Vazquez (see \cite{Vaz}), $u_{n, \epsilon}$ is positive everywhere.
Hence, setting $\sigma_{n, \epsilon}=|\nabla u_{n,\epsilon}|^{\epsilon-1} \nabla u_{n,\epsilon}$, we have shown that  $u_{n,\epsilon}\in   \mathcal C^{1,\alpha}(\bar \Omega)\cap W^{1,1+\epsilon}_0 (\Omega)$ is a positive solution of (\ref{1Lap:eq:P_{n,epsilon}}).

\begin{lemma}\label{1Lap:lemma:1.1}
The problem {\rm (\ref{1Lap:eq:P_{n,epsilon}})} has a unique positive solution
\end{lemma}

\begin{proof}[Proof of Lemma \ref{1Lap:lemma:1.1}] 
Let $u$ and $v$ be two positive solutions of (\ref{1Lap:eq:P_{n,epsilon}}). Then we have:
\begin{equation}\label{1Lap:eq:2.8}
-{\rm div}\left[\sigma_{n, \epsilon}(u)-\sigma_{n, \epsilon}(v)\right]+
2n\left[\alpha(u)-\alpha(v)\right]u^{\epsilon}+2n \  \alpha(v)\left(u^{\epsilon}-v^{\epsilon}\right)=0,
\end{equation}
where $\alpha(u)={\ds \int_\Omega} u^{1+\epsilon} -1.$

\bigskip

{\bf Case 1:}
$\|u\|_{1+\epsilon} \geq \|v\|_{1+\epsilon}.$\\   
Let us multiply  (\ref{1Lap:eq:2.8}) by $(u-v)^+$ then integrate.
It is clear that 
$$2n\left[\alpha(u)-\alpha(v)\right]\int_\Omega u^{\epsilon}(u-v)^+ \geq 0 .$$
So we get that:
\begin{equation}\label{1lap:eq:2.11}
\int_\Omega[\sigma_{n, \epsilon}(u)-\sigma_{n, \epsilon}(v)] \cdot \nabla(u-v)^+ +2n \  \alpha(v)\int_\Omega\left(u^{\epsilon}-v^{\epsilon}\right)(u-v)^+ \leq0 .
\end{equation}
We know that 
\begin{equation}\label{1Lap:eq:2.9}
\int_\Omega [\sigma_{n, \epsilon}(u)-\sigma_{n, \epsilon}(v)] \cdot \nabla(u-v)\geq 0.
\end{equation}
On the other hand it is clear that
\begin{equation}\label{1Lap:eq:2.10} 
\int_\Omega \left(u^{\epsilon}-v^{\epsilon}\right)(u-v)\geq 0.
\end{equation}
So, we can conclude that:
\begin{equation}\label{1lap:eq:2.12}
\int_\Omega[\sigma_{n, \epsilon}(u)-\sigma_{n, \epsilon}(v)] \cdot \nabla(u-v)^++2n \  \alpha(v)\int_\Omega \left(u^{\epsilon}-v^{\epsilon}\right)(u-v)^+ \geq0 .
\end{equation}
So from (\ref{1lap:eq:2.11}) and  (\ref{1lap:eq:2.12}), we obtain that
$$\int_\Omega[\sigma_{n, \epsilon}(u)-\sigma_{n, \epsilon}(v)] \cdot \nabla(u-v)^++2n \  \alpha(v)\int_\Omega\left(u^{\epsilon}-v^{\epsilon}\right)(u-v)^+ =0 .$$
Then ${\ds \int_\Omega} \left(u^{\epsilon}-v^{\epsilon}\right)(u-v)^+=0$, which implies $(u-v)^+=0$, i.e. $u\leq v$.
Using $\|u\|_{1+\epsilon} \geq\|v\|_{1+\epsilon}$, one finally gets $u=v$ a.e.

\bigskip

{\bf Case 2:}
$\|u\|_{1+\epsilon} \leq\|v\|_{1+\epsilon}$.\\
We use the same arguments as in the Case 1, just replacing $(u-v)^+$ by $(v-u)^+$.   
\end{proof}
Thus, we have proved the existence and uniqueness of a positive solution to the problem (\ref{1Lap:eq:2.2}).

\paragraph{Step 2:}
We prove here that ${\ds \lim_{\epsilon \rightarrow 0}}\lambda_{1+\epsilon, n}=\lambda_{1, n}$.
\begin{proposition}\label{1Lap:prop:1.2}
For every $n \in \N^*$, we have:
$$\limsup_{\epsilon \rightarrow 0} \ \lambda_{1+\epsilon,  n} \leq \lambda_{1, n}$$
\end{proposition}

\begin{proof}[Proof of Proposition \ref{1Lap:prop:1.2}]
Let $\delta>0$ be given and $\phi \in \mathcal D(\Omega)$ such that
$$I_{1, n}(\phi)=\int_\Omega |\nabla \phi | +  n \left( \int_\Omega| \phi|-1 \right)^2 \leq \lambda_{1, n}+\delta.$$
But ${\ds \lim_{\epsilon \rightarrow 0}} I_{1+\epsilon, n}(\phi)=I_{1, n}(\phi)$,
hence,
$$\limsup_{\epsilon \rightarrow 0} \lambda_{1+\epsilon,  n} \leq \lambda_{1, n}+\delta.$$
$\delta$ being arbitrary, we get ${\ds \limsup_{\epsilon \rightarrow 0}} 
\lambda_{1+\epsilon,  n} \leq \lambda_{1, n}.$ 
\end{proof}
Let now $u_{n,\epsilon}$ be the positive solution of the minimizing problem (\ref{1Lap:eq:2.2}). 
Using Poincar\'e's and H\"older's inequalities, we get
\begin{align*}
\int_\Omega u_{n,\epsilon} dx 
\leq  C \int_\Omega |\nabla u_{n,\epsilon}| dx  
&\leq C\left(\int_\Omega 
|\nabla u_{n,\epsilon}|^{1+\epsilon} dx\right)^{\frac{1}{1+\epsilon}} 
|\Omega|^{\frac{\epsilon}{1+\epsilon}}\\ 
&\leq  C^\prime \lambda_{1+\epsilon,  n}.
\end{align*}
These inequalities show that $(u_{n,\epsilon})_{\epsilon>0}$ and $(\nabla u_{n,\epsilon})_{\epsilon>0}$
are both bounded in $L^1(\Omega)$.
This means that $(u_{n,\epsilon})_{\epsilon>0}$ is bounded in $BV(\Omega)$. 
We denote by $u_n$ the limit of some subsequence in $BV$ for the weak topology.  

In step 4 we shall precise this limit. In particular we shall  obtain $u_n$ as  the restriction to $\Omega$ of  some limit of extended functions $u_{n,\epsilon}$ by zero  outside of $\Omega$.

\paragraph{Step 3:} we obtain $\sigma_n=\frac{``\nabla u_n"}{| \nabla u_n|}$
as the weak limit of $\sigma_{n,\epsilon}=| \nabla u_{n,\epsilon}|^{\epsilon-1} \nabla u_{n,\epsilon}$.\\

Let $\sigma_{n,\epsilon}=| \nabla u_{n,\epsilon}|^{\epsilon-1} \nabla u_{n,\epsilon}$, 
one sees that $\sigma_{n,\epsilon}$ is uniformly bounded in $L^{\frac{1+\epsilon}{\epsilon}}(\Omega)$.
Let us prove that $\sigma_{n, \epsilon}$ is uniformly bounded in every
$L^{q}(\Omega)$,  $\text{for all} \  q< \infty$.
Indeed, let $q>1$ be given and let $\epsilon$ be such that $q< \frac{1+\epsilon}{\epsilon}$.
Then
$$\left(\int_\Omega |\sigma_{n,\epsilon}|^ q\right)^{\frac{1}{q}} 
\leq \left(\int_\Omega |\sigma_{n,\epsilon}|^{\frac{1+\epsilon}{\epsilon}}\right)^{\frac{\epsilon}{1+\epsilon}} |\Omega|^{\frac{1+\epsilon(1-q)}{(1+\epsilon)q}}\leq C.$$
Then we may extract from it a subsequence, still denoted by $\sigma_{n,\epsilon}$, such that
$\sigma_{n,\epsilon}$ tends to some $\sigma_n$ weakly in $L^{q}(\Omega)$,  $\text{for all} \  q< \infty$ and $\sigma_{n,\epsilon}$ tends to $\sigma_n$ a.e., when $\epsilon$ tends to $0$. 

We observe now that $\|\sigma_n\|_{\infty} \leq 1$. 
For that aim, let $\eta$ be in $\mathcal D(\Omega, \R^N)$. Then
\begin{align*}
\left| \int_\Omega \sigma_n \cdot \eta \right| 
\leq {\ds \liminf_{\epsilon\rightarrow 0}}\left|\int_\Omega \sigma_{n,\epsilon} \cdot \eta\right| 
& \leq {\ds \liminf_{\epsilon \rightarrow 0}} \int_\Omega|\nabla u_{n,\epsilon}|^{\epsilon}|\eta|\\
&\leq {\ds \liminf_{\epsilon \rightarrow 0}} \left(\int_\Omega| \nabla u_{n,\epsilon}|^{1+\epsilon} \right)^{\frac{\epsilon}{1+\epsilon}} \left(\int_\Omega|\eta|^{1+\epsilon} \right)^{\frac{1}{1+\epsilon}}\\
&\leq {\ds \liminf_{\epsilon \rightarrow 0}} \left(\lambda_{1+\epsilon,  n}\right)^
{\frac{\epsilon}{1+\epsilon}} \left(\int_\Omega|\eta|^{1+\epsilon} \right)^
{\frac{1}{1+\epsilon}}\\
&\leq \int_\Omega |\eta|.
\end{align*}
This implies that  $\|\sigma_n\|_{\infty} \leq 1$.

Let us now observe that $u_{n, \epsilon}^\epsilon$ is uniformly  bounded in every $L^q(\Omega),  q< \infty$. Indeed, let $q$ be given and  let $\epsilon$ be small enough, such that $q< \frac{1+\epsilon}{\epsilon}$, then 
$$\left(\int_\Omega |u_{n,\epsilon}^{\epsilon}|^ q \right)^{\frac{1}{q}}
\leq \left(\int_\Omega |u_{n,\epsilon}|^{1+\epsilon}\right)^{\frac{\epsilon}{1+\epsilon}}
|\Omega|^{\frac{1+\epsilon(1-q)}{q(1+\epsilon)}} \leq C.$$
Then $w_{n,\epsilon} = u_{n,\epsilon}^\epsilon$ converges  weakly, in every $L^q(\Omega), \ q<\infty$, up to a subsequence, to some $w_n$,  when $\epsilon$ tends to $0$. 

Let us prove that $0\leq w_n\leq 1$ and $(w_n-1) u_n = 0$. 
For the first assertion, let $\eta \in \mathcal D(\Omega)$,
\begin{align*}
\left|\int_\Omega w_n \cdot \eta \right| 
&\leq \left(\int_\Omega|w_n|^{\frac{1+\epsilon}{\epsilon}} \right)^{\frac{\epsilon}{1+\epsilon}} \left(\int_\Omega|\eta|^{1+\epsilon} \right)^{\frac{1}{1+\epsilon}}\\
&\leq \liminf_{\epsilon \rightarrow 0} \left( \lambda_{1+\epsilon,  n}\right)^
{\frac{\epsilon}{1+\epsilon}} \left(\int_\Omega|\eta|^{1+\epsilon} \right)^
{\frac{1}{1+\epsilon}}\\&\leq \int_\Omega |\eta|.
\end{align*}
Hence $0\leq w_n\leq 1, \quad \forall n \in \N^*.$

To prove that $(w_n-1) u_n = 0$, let us observe that 
$u_{n, \epsilon} \longrightarrow u_n$ in $L^k(\Omega)$ strongly for all 
$k< \frac{N}{ N-1}$ and $w_{n,\epsilon} \longrightarrow w_n$ in $L^{N+1}(\Omega)$ weakly, therefore 
$$\int_\Omega  w_{n,\epsilon} u_{n,\epsilon} \longrightarrow \int_\Omega w_n u_n \quad \text{when} \ \epsilon  \rightarrow 0$$
Finally,
$$\int_\Omega w_n u_n 
= \lim_{\epsilon \rightarrow 0}\int_\Omega  w_{n,\epsilon} u_{n,\epsilon}
=\lim_{\epsilon \rightarrow 0}\int_\Omega u_{n,\epsilon}^{1+\epsilon}
=\int_\Omega u_n.$$
Using the fact that $0\leq w_n\leq 1$, one gets the result. 

Passing to the limit in (\ref{1Lap:eq:2.5}), one gets:
\begin{equation}\label{1Lap:eq:2.13}
-{\rm div} \sigma_n +2n \left( \int_ \Omega u_n -1 \right) w_n=0 .
\end{equation}

\paragraph{Step 4:} Extension of $u_{n,\epsilon}$ outside $\Omega$ and convergence towards a solution of (\ref{1Lap:eq:P_{n,epsilon}}).\\

Let $\widetilde u_{n,\epsilon}$ be the extension of $u_{n,\epsilon}$ by $0$ in $\R^N\setminus \overline{\Omega}$. 
Since $u_{n,\epsilon}=0$ on $\partial \Omega$, then $\widetilde u_{n,\epsilon} \in W^{1,1+\epsilon}(\R^N)$ 
and $(\widetilde u_{n,\epsilon})$ is bounded in  $BV(\R^N)$.
Then one may extract from it a subsequence, still denoted $(\widetilde u_{n,\epsilon})$, such that
$$\widetilde u_{n,\epsilon} \longrightarrow v_n  \quad \text{in} \ L^k (\R^N), \quad \forall \ k< \frac{N}{N-1} \quad \text{when} \ \epsilon \longrightarrow 0,$$
with $v_n=0$ outside of $\overline{\Omega}$
and $$\nabla \widetilde u_{n,\epsilon} \rightharpoonup  \nabla v_n  \quad \text{weakly \ in} \ {\cal M}^1(\R^N) \quad \text{when} \ \epsilon \longrightarrow 0,$$
We denote by $u_n$ the restriction of $v_n$ to  $\Omega$. 
 We use in the above  some  limit $\sigma_n$ of  
$\sigma_{n,\epsilon}=|\nabla u_{n,\epsilon}|^{\epsilon-1} \nabla u_{n,\epsilon}$ obtained in the third step.

Multiplying the equation (\ref{1Lap:eq:2.5}) by $\widetilde u_{n,\epsilon} \phi$, where $\phi \in \ \mathcal D(\R^N)$, and integrating by parts, one obtains:
$$\int_{{\Omega}} \sigma_{n,\epsilon} \cdot \nabla (\widetilde u_{n,\epsilon} \phi)+2n \left( \int_ \Omega \widetilde u_{n,\epsilon}^{1+\epsilon} -1\right)\int_ \Omega \widetilde u_{n,\epsilon}^{1+\epsilon} \phi =0 \ , $$
or equivalently 
\begin{equation}\label{1Lap:eq:equivalence}
\int_{\R^N} |\nabla (\widetilde u_{n,\epsilon})|^{1+\epsilon} \phi+\int_{\Omega}\sigma_{n,\epsilon}u_{n,\epsilon} \cdot \nabla \phi+2n \left(\int_{\R^N}\widetilde u_{n,\epsilon}^{1+\epsilon} -1\right) \int_{\R^N}\widetilde u_{n,\epsilon}^{1+\epsilon} \phi =0 .
\end{equation}
Since $\sigma_{n,\epsilon} \rightharpoonup \sigma_n$ in  $L^q (\Omega)$ for all $q< \infty$, in particular for any $\alpha >0$, $\sigma_{n,\epsilon}$ tends weakly towards $\sigma_n$ in  $L^{N+\alpha} (\Omega)$.
Since $\widetilde u_{n,\epsilon}$ tends strongly towards $v_n$ in $L^k(\Omega), \ k< \frac{N}{N-1}$, one obtains that:
$$\int_{\Omega}\sigma_{n,\epsilon}u_{n,\epsilon} \cdot \nabla \phi \longrightarrow \int_{\Omega}\sigma_nu_n \cdot \nabla \phi,  \ \quad \text{when} \ \epsilon\rightarrow 0 . $$
By passing to the limit in the equation (\ref{1Lap:eq:equivalence})
and defining, up to extracting a subsequence, the measure $\mu$ on $\R^N$ by: $\lim_{\epsilon \rightarrow 0}|\nabla (\widetilde u_{n,\epsilon})|^{1+\epsilon}=\mu$, one obtains:
\begin{equation}\label{1Lap:eq:2.14}
\langle \mu, \phi \rangle+ \int_{\Omega}\sigma_nu_n \cdot \nabla \phi+2n \left(\int_{\R^N}v_n-1\right) \int_{\R^N}v_n \phi =0 .
\end{equation}
On the other hand, multiplying equation (\ref{1Lap:eq:2.13}) by $v_{n} \phi$ where $\phi \in  \mathcal D(\R^N)$, one gets
\begin{equation}\label{1Lap:eq:2.15}
\int_{\Omega \cup \partial \Omega}\sigma_n \cdot (\nabla v_n) \phi
+\int_{\Omega}\sigma_n u_n \cdot \nabla \phi
+2n \left(\int_{\Omega}u_n -1\right)  \int_{\Omega}  u_n \phi=0.
\end{equation}
Substracting (\ref{1Lap:eq:2.15}) from (\ref{1Lap:eq:2.14}), one gets 
\begin{equation}\label{1Lap:eq:2.16}
\mu=\sigma_n \cdot \nabla v_n \quad \text{in} \ \Omega \cup \partial \Omega.
\end{equation}
This implies in particular, according to Proposition \ref{1Lap:prop:1.1}, that
$$|\mu|\leq |\nabla v_n| \quad \text{in} \ \Omega \cup \partial \Omega,$$
and
$$\int_{\R^N} |\nabla (\widetilde u_{n,\epsilon})|^{1+\epsilon} \longrightarrow \int_{\R^N} |\nabla v_n| \quad \text{when} \ \epsilon  \rightarrow 0.$$
Finally recalling that according to proposition \ref{1Lap:prop:1.1}  , one has $\nabla  v_n \cdot \sigma_n\leq |\nabla v_n|$ on $\Omega\cup\partial\Omega$ one derives that 
$$|\nabla v_n|= \sigma_n \cdot \nabla v_n \quad \text{in} \ \Omega \cup \partial \Omega.$$
Recall that  from Proposition \ref{1Lap:prop:1.1}
$$\nabla v_n = \nabla u_n \chi_{\Omega} -  u_n  \ \delta_{\partial \Omega} \overrightarrow {n} , $$
$$\sigma_n \cdot \nabla v_n=\sigma_n \cdot \nabla u_n \chi_{\Omega} -\sigma_n \cdot \overrightarrow {n} u_n \delta_{\partial \Omega} , $$
we have obtained 
$$\begin{cases}
\sigma_n \cdot \nabla u_n=|\nabla u_n|  \ \quad \text{in} \  \Omega , \\
\sigma_n \cdot \overrightarrow {n} u_n = - u_n \ \quad \text{on} \ \partial \Omega . 
\end{cases}$$
Then $u_n$ is a nonnegative solution of (\ref{1Lap:eq:P_n}). Moreover, the convergence of $|\nabla \widetilde u_{n,\epsilon}|$ is tight on $\overline{\Omega}$, i.e.
$$\int_{\Omega}|\nabla  u_{n,\epsilon}| \longrightarrow\int_{\Omega}|\nabla u_n|+\int_{\partial \Omega} u_n, \ \quad \text{when} \ \epsilon\rightarrow 0 . $$
Indeed, one has ${\ds \int_\Omega} |\nabla u_{n,\epsilon}|^{1+\epsilon}\longrightarrow {\ds \int_\Omega} |\nabla u_n|+ {\ds \int_{\partial \Omega }}u_n \quad \text{when} \ \epsilon \longrightarrow 0$. Using the lower semicontinuity for the extension $u_{n,\epsilon}$ and H\"older's inequality, we get
\begin{align*}
{\ds \lim_{\epsilon \rightarrow 0}}\int_\Omega |\nabla u_{n,\epsilon}|^{1+\epsilon}
=\int_\Omega |\nabla u_n|+\int_{\partial \Omega }u_n
&\leq \liminf_{\epsilon\rightarrow 0} \int_{\Omega} |\nabla u_{n,\epsilon}|\\
&\leq \liminf_{\epsilon\rightarrow 0} \left(\int_{\Omega} |\nabla u_{n,\epsilon}|^{1+\epsilon}\right)^{\frac{1}{1+\epsilon}} |\Omega |^{\frac{\epsilon}{1+\epsilon}}\\
&={\ds \lim_{\epsilon \rightarrow 0}}\int_\Omega |\nabla u_{n,\epsilon}|^{1+\epsilon}
\end{align*}
The result is proved.

\paragraph{Step 5:} The convergence of $\lambda_{1+\epsilon, n}$ towards $\lambda_1$\\

In this step we explicit the relation between the values $\lambda_{1+\epsilon, n}$ when $n$ is large,  and the first eingenvalue $\lambda_1$ defined in the first part.

\begin{theorem}\label{1Lap:thm:1.2}
Let $u_n$ be  a nonnegative solution of  \ref{1Lap:eq:2.5},  then, up to a subsequence, as $n\rightarrow \infty$, $(u_n)$ converges to $u \in BV(\Omega), \ u \geq0, \ u \not \equiv 0,$ which realizes the minimum defined in {\rm (\ref{1Lap:eq:I_BV})}.
Moreover
$$\lim_{n \rightarrow \infty} \lambda_{1,  n} = \lambda_1 . $$ 
\end{theorem}

\begin{proof}[Proof of the Theorem \ref{1Lap:thm:1.2}]
 For $\lambda_{1,  n}$ and $\lambda_1$ defined as above, it is clear that we have:
 \begin{equation}\label{1Lap:eq:3.3}
 \limsup_{n \rightarrow \infty} \lambda_{1,  n} \leq \lambda_1.
 \end{equation}
Let $(u_n)_n$ be a sequence of positive solutions  of the relaxed problem defined in (\ref{1Lap:eq:lambda_BV}). We begin to prove that $(u_n)_n$ is bounded in $BV(\Omega)$. For that aim let us note that by (\ref{1Lap:eq:3.3}), one gets that  $n\left({\ds \int}_\Omega u_n -1 \right)^2$ is bounded by $\lambda_1$, which implies that $\lim_{n \rightarrow \infty} \left({\ds \int}_\Omega u_n-1\right)^2 = 0$.
Then
$$\lim_{n \rightarrow \infty} \int_\Omega u_n = 1 , $$
Hence, $(u_n)_n$ is bounded in $L^1(\Omega)$.\\
Using once more (\ref{1Lap:eq:3.3}), one can conclude that $(u_n)_n$ is bounded in $BV(\Omega)$.
\noindent Then, the extension of each $u_n$ by zero outside of $\overline{\Omega}$ is bounded in $BV(\R^N)$. One can then extract from it a subsequence, still denoted $u_n$, such that 
$$u_n \rightharpoonup u  \quad \text{weakly \ in} \ BV(\R^N) \quad \text{when} \ n  \rightarrow \infty,$$
obviously $u=0$ outside of $\overline{\Omega}$ and $u > 0$ in $\Omega$.  By the compactness of the Sobolev embedding from  $BV(\Omega)$ into $L^1(\Omega)$, one has $\|u\|_{L^1(\Omega)}=1$. 
Using the lower semi continuity of th total variation $\int_{\R^N} |\nabla u|$ with respect to the weak topology, one has (since $u_n \rightarrow u$ in $L^1(\Omega)$)
\begin{align*}
\lambda_1 \leq \int_{\R^N} |\nabla u| 
&\leq  \int_{\R^N}|\nabla u |+ n\left(\int_{\R^N}u -1\right)^2 \\
&\leq \liminf_{n \rightarrow \infty} \left[\int_{\R^N}|\nabla u_n|+ n\left(\int_{\R^N} u_n -1\right)^2 \right] \\
&\leq \limsup_{n \rightarrow \infty} \lambda_{1,  n} \leq \lambda_1.
\end{align*}
Then one gets that 
$$\lim_{n \rightarrow \infty} \lambda_{1,  n} =\lambda_1 . $$
Since $u=0$ outside of $\overline{\Omega}$, one has  on $\R^N$ $\nabla u =\nabla u \chi_{\Omega} - u   \overrightarrow {n} \ \delta_{\partial \Omega}$ and then $$\int_{\R^N}|\nabla u |=\int_{\Omega}|\nabla u|+\int_{\partial \Omega}u . $$
Moreover, one obtains that:
$$\lim_{n \rightarrow \infty} n \left( \int_\Omega u_n-1\right)^2=0 , $$
and 
$$\lim_{n \rightarrow \infty} \int_\Omega|\nabla u_n|= \int_\Omega|\nabla u|+\int_{\partial \Omega}|u| . $$
Then, we get the tight convergence of $u_n$ to $u$ in $BV(\overline{\Omega})$. 

Let us observe that ${\rm sign}^+ (u_n)$ converges to some $w$, $0\leq w\leq 1$ in every $L^q(\Omega), \forall \ q< \infty$. Using the convergence of $u_n$ to $u$ in $L^q(\Omega), \ \forall \ q< \frac{N}{N-1}$, one gets 
$$\int_\Omega u_n = \int_\Omega u_n {\rm sign}^+ (u_n) \longrightarrow  \int_\Omega u=1 \quad \text{when} \ \ n \rightarrow \infty.$$
As a consequence
$$-2n \left( \int_\Omega u_n -1 \right)\int_\Omega u_n  \longrightarrow \lambda_1 \quad \text{when} \ \ n \rightarrow \infty,$$
and then also
$$-2n \left( \int_\Omega u_n -1 \right) \longrightarrow \lambda_1 \quad \text{when} \ \ n \rightarrow \infty.$$
This ends the proof of the main result.
\end{proof}

The author thanks the referee for its remarks and advices which permit to improve this paper.

\index{}

\end{document}